\documentclass{amsart}
\usepackage{amssymb}
\begin{document}
\newtheorem{theorem}{Theorem}[section]
\newtheorem{lemma}[theorem]{Lemma}
\newtheorem{remark}[theorem]{Remark}
\newtheorem{definition}[theorem]{Definition}
\newtheorem{corollary}[theorem]{Corollary}
\newtheorem{example}[theorem]{Example}
\newtheorem{conjecture}[theorem]{Conjecture}
\font\pbglie=eufm10
\def\qedbox{\hbox{$\rlap{$\sqcap$}\sqcup$}}
\def\RRR{{\text{\pbglie R}}}
\def\BB{{\mathcal{B}}}
\def\Tr{\operatorname{Tr}}
\def\id{\operatorname{Id}}
\def\ffrac#1#2{{\textstyle\frac{#1}{#2}}}
\def\Span{\operatorname{Span}}
\def\Range{\operatorname{Range}}
\def\Rank{\operatorname{Rank}}
\def\text#1{{\hbox{#1}}}
\def\operatorname#1{{\rm#1\,}}
\makeatletter
 \renewcommand{\theequation}{%
 \thesection.\alph{equation}}
 \@addtoreset{equation}{section}
 \makeatother
\author{N. Blazic, P. Gilkey, S. Nikcevic, and U. Simon}
\title{The spectral geometry of the Weyl conformal tensor}
\begin{address}{NB: Faculty of Mathematics, University of Belgrade\\ Studentski Trg. 16, P.P. 550,
11000 Beograd, Srbija i Crna Gora.\\ Email: {\it blazicn@matf.bg.ac.yu}}\end{address}
\begin{address}{PG:Mathematics Department, University of Oregon\\
Eugene Or 97403 USA.\newline Email: {\it gilkey@darkwing.uoregon.edu}}
\end{address}
\begin{address}{SN: Mathematical Institute, SANU\\
Knez Mihailova 35, p.p. 367,
11001 Beograd,
Srbija i Crna Gora.\\
Email: {\it stanan@mi.sanu.ac.yu}}\end{address}
\begin{address}
{US: Fachbereich Mathematik, Technische Universit\"at
Berlin\\ Strasse des 17. Juni 135, D-10623 Berlin, Deutschland\\
Email: {\it simon@sfb288.math.tu-berlin.de}}\end{address}
\subjclass[2000]{Primary 53C50\\\hphantom{....}
{\it Key words:} conformally flat, Ivanov-Petrova manifold, Jacobi operator, Osserman manifold, skew-symmetric curvature
operator, Weyl conformal tensor.}
\begin{abstract}{We study when the Jacobi operator
associated to the Weyl conformal curvature tensor has constant eigenvalues on the bundle of unit spacelike or timelike tangent vectors.
This leads to questions in the conformal geometry of pseudo-Riemannian manifolds which generalize the Osserman conjecture
to this setting. We also study similar questions related to the skew-symmetric curvature operator defined by the Weyl conformal
curvature tensor.}\end{abstract}
\maketitle
 
\section{Introduction}\label{Sect-1}

\subsection{Algebraic curvature tensors} We work in a purely algebraic context for the moment. Consider a triple
$\mathcal{V}:=(V,g,A)$ where $g$ is a non-degenerate inner product of signature $(p,q)$ on a finite dimensional real vector space $V$
of dimension $m:=p+q\ge3$ and where
$A\in\otimes^4V^*$ is an {\it algebraic curvature tensor} on $V$; i.e. $A$ has the usual symmetries of the Riemann
curvature tensor:
\begin{eqnarray*}
&&A(x,y,z,w)=A(z,w,x,y)=-A(y,x,z,w),\quad\text{and}\\
&&A(x,y,z,w)+A(y,z,x,w)+A(z,x,y,w)=0\,.
\end{eqnarray*}
We say that $\mathcal{V}$ is {\it Riemannian} if $p=0$ and {\it Lorentzian} if $p=1$.

If $\phi$ is a $g$ self-adjoint endomorphism of $V$, then we set
\begin{equation}\label{eqn-1.a}
A_\phi(x,y,z,w):=g(\phi x,w)g(\phi y,z)-g(\phi x,z)g(\phi y,w)\,.
\end{equation}
Fiedler \cite{F01} showed these algebraic curvature tensors span the space $\mathcal{A}(V)$ of all algebraic curvature tensors. These
tensors will play a crucial role in our subsequent discussion; setting $\phi=\id$ yields the algebraic curvature tensor 
of constant
sectional curvature $+1$.

\subsection{The Weyl conformal curvature tensor} There is
a natural representation of the orthogonal group
$O(V,g)$ on $\mathcal{A}(V)$ defined by pull-back; if $A\in\mathcal{A}(V)$ and $g\in O(V,g)$, the pull-back $g^*A\in\mathcal{A}(V)$ is given
by
$$(g^*A)(x,y,z,w)=A(gx,gy,gz,gw)\,.$$
This representation is not irreducible but decomposes as the direct sum of 3 irreducible representations which we can describe as follows.
Let $g_{ij}:=g(e_i,e_j)$ and let $g^{ij}$ be the inverse matrix relative to some basis $\{e_i\}$ for $V$. The associated Ricci tensor
$\rho_A$ and scalar curvature $\tau_A$ are then defined by contracting indices:
$$\rho_A(x,y):=\textstyle\sum_{ij}g^{ij}A(x,e_i,e_j,y)\quad\text{and}\quad\tau_A:=\textstyle\sum_{ij}g^{ij}\rho_A(e_i,e_j)\,.$$
The associated maps $\sigma_\rho:A\rightarrow\rho_A\in S^2(V^*)$ and $\sigma_\tau:A\rightarrow\tau_A\in\mathbb{R}$ are
$O(V,g)$ equivariant. The space
$\mathcal{W}(V,g):=\ker(\sigma_\rho)$ of algebraic Weyl tensors is an irreducible
representation space for
$O(V,g)$ and we have:
$$\mathcal{A}(V)=\mathcal{W}(V,g)\oplus S^2(V^*)$$
as an $O(V,g)$ representation space. The further decomposition of $S^2(V^*)$ as the direct sum of the trace free
tensors and the scalar multiples of the identity then completes the decomposition of $\mathcal{A}(V)$ as a direct sum of irreducible $O(V,g)$
modules.
Let $\pi_{\mathcal{W}}$ be orthogonal projection from $\otimes^4V^*$ to $\mathcal{W}(V,g)$:
\begin{eqnarray}\nonumber
\pi_{\mathcal{W}}(A)(x,y,z,w)&=&
A(x,y,z,w)-\ffrac1{m-2}\{\rho_A(x,w)g(y,z)+g(x,w)\rho_A(y,z)\}\nonumber\\
&+&\ffrac1{m-2}\{\rho_A(x,z)g(y,w)+g(x,z)\rho_A(y,w)\}\label{eqn-1.b}\\
&+&\ffrac1{(m-1)(m-2)}\tau_A\{g(x,w)g(y,z)-g(x,z)g(y,w)\}\,.\nonumber
\end{eqnarray}

\subsection{The Jacobi operator} If $A$ is an algebraic curvature tensor, then the {\it Jacobi operator} $J_A$ is a $g$ self-adjoint map of
$V$ characterized by the property:
$$g(J_A(x)y,z)=A(y,x,x,z)\,.$$
For example, if $A=A_\phi$ is given by Equation (\ref{eqn-1.a}), then
\begin{equation}\label{eqn-1.c}
J_{A_\phi}(x)y=g(\phi x,x)\phi y-g(\phi x,y)\phi x\,.
\end{equation}

It is clear that $\rho_A(x,x)=\Tr\{J_A(x)\}$ for any $A$;
in particular
\begin{equation}\label{eqn-1.d}
\Tr\{J_W(x)\}=0\quad\text{for any}\quad x\in V\quad\text{if}\quad W\in\mathcal{W}(V,g)\,.
\end{equation}
The pseudo-spheres of unit spacelike ($+$) and unit timelike ($-$) vectors in $V$ are
$$S^\pm(\mathcal{V}):=\{v\in V:g(v,v)=\pm1\}\,.$$
We say that $\mathcal{V}$ is {\it spacelike} (resp. {\it timelike}) {\it Jordan Osserman} if the Jordan normal form of
$J_A$ is constant on $S^+(\mathcal{V})$ (resp. on $S^-(\mathcal{V})$).
If
$\mathcal{V}$ is Riemannian, then the Jordan normal form is determined by the eigenvalue structure and, as every non-zero vector is
spacelike, we shall drop the qualifiers `spacelike' and `Jordan' in the interests of notational simplicity. Note that the eigenvalue
structure does not determine the Jordan normal form in the higher signature context.

\subsection{The skew-symmetric curvature operator} Let
$\{e_1,e_2\}$ be an orthonormal basis for an oriented spacelike (resp. timelike)
 $2$ plane $\pi$ of $V$. One then defines the {\it skew-symmetric} curvature operator $A(\pi)$ by the identity:
$$g(A(\pi)x,y):=A(e_1,e_2,x,y)\,.$$
This $g$ skew-symmetric endomorphism of $V$ is independent of the particular oriented orthonormal basis for $\pi$ which is chosen. One
says
$\mathcal{V}$ is {\it spacelike} (resp. {\it timelike}) {\it Jordan Ivanov-Petrova} if the Jordan normal form of $A(\pi)$ is
constant on the Grassmannian of oriented spacelike (resp. timelike) $2$ planes in $V$. 

\subsection{The geometric setting} Let $R$ be the Riemann curvature tensor of a pseudo-Riemannian manifold $(M,g)$ of signature
$(p,q)$ and dimension $m:=p+q\ge3$. Let $\mathcal{R}_P:=(T_PM,g_P,R_P)$ be the triple determined by the tangent
bundle of $M$ at a point $P$ of $M$, the pseudo-Riemannian metric $g_P$, and the curvature tensor $R_P$. 

We say that $(M,g)$ is {\it pointwise spacelike} (resp. {\it timelike}) {\it Jordan Osserman} if $\mathcal{R}_P$ is
spacelike (resp. timelike) Jordan Osserman for every point $P$ of $M$; the Jordan normal form of $J_R$ is allowed to vary with the point $P$
of $M$. We say
$(M,g)$ is globally spacelike (resp. timelike) Jordan Osserman if the Jordan normal form of $J_R$ on the appropriate pseudo-sphere bundle is
independent of $P$. It is known that any global Riemannian ($p=0$) Osserman manifold is locally isometric to a rank $1$ symmetric space if
$m\ne8,16$ \cite{Ch88,Ni03,refNik2} and that any local Lorentzian ($p=1$) Jordan Osserman manifold has constant sectional
curvature \cite{BBG97,GKV97}. In the higher signature setting, there exist spacelike and timelike Jordan Osserman manifolds which are not
locally homogeneous \cite{BBZ01,GVV98}. There is a vast literature on the subject and we shall content ourselves by refering to
\cite{GKV02} for further details.

We say that $(M,g)$ is {\it pointwise spacelike} (resp. {\it timelike}) {\it Jordan Ivanov-Petrova} if $\mathcal{R}_P$ is spacelike
(resp. timelike) Jordan Ivanov-Petrova for every point $P$ of $M$; again, the Jordan normal form of $J_R$ is allowed to very with the point
$P$ of $M$. These manifolds have been classified in the Riemannian setting if $m\ne 3,7$ \cite{Gi99,GLS99,IvPe98}, in the Lorentzian
setting if $m\ge11$ and if $\{m,m+1\}$ are not powers of $2$ \cite{GiZa02a,GiZa02b}, and in the higher signature setting if $q\ge11$, if
$p\le\frac{q-6}4$, if
$\{q,q+1,...,q+p\}$ does not contain a power of $2$, and if $R(\pi)$ is not nilpotent \cite{St03}. We refer to \cite{refGi02} for further
details concerning spacelike and timelike Jordan Ivanov-Petrova manifolds.

\subsection{Conformal geometry} Let $P$ be a point of a pseudo-Riemannian manifold $(M,g)$. Let
$\mathcal{W}_P:=(T_PM,g_P,W_P)$ where $W_P:=\pi_{\mathcal{W}}R_P$ is the associated Weyl conformal curvature tensor on $T_PM$. We say that
$(M,g)$ is {\it conformally spacelike} (resp. {\it timelike}) {\it Jordan Osserman} if $\mathcal{W}_P$ is spacelike (resp. timelike) Jordan
Osserman for every point $P$ of $M$. Similarly, we say that $(M,g)$ is {\it conformally spacelike} (resp. {\it timelike}) {\it Jordan
Ivanov-Petrova} if $\mathcal{W}_P$ is spacelike (resp. timelike) Jordan Ivanova-Petrova for every point $P$ of
$M$. In both settings, the Jordan normal form is permitted to vary with the point $P$ of $M$; the technical distinction between `global' and
`pointwise' plays no role in this setting.

Recall that two metrics $g_1$ and $g_2$ are said to be {\it conformally equivalent} if there is a positive scaling function
$\alpha\in C^\infty(M)$ so that $g_1=\alpha g_2$. We let $[g]$ be the set of all pseudo-Riemannian metrics on $M$ which are conformally
equivalent to $g$.

\begin{theorem}\label{thm-1.1}
Let $g_1\in[g_2]$. Then:
\begin{enumerate}
\item $(M,g_1)$ is conformally spacelike (resp. timelike)
Jordan Osserman if and only if $(M,g_2)$ is conformally spacelike (resp. timelike) Jordan Osserman.
\item $(M,g_1)$ is conformally spacelike (resp. timelike)
Jordan Ivanov-Petrova if and only if $(M,g_2)$ is conformally spacelike (resp. timelike) Jordan Ivanov-Petrova.
\end{enumerate}
\end{theorem}

\begin{proof} As $g_1=\alpha g_2$, one has $W_{g_1}=\alpha W_{g_2}$; the Weyl conformal curvature tensor simply rescales. Let $x\in T_PM$
be a
$g_2$ spacelike or timelike unit vector. Let
$$\tilde x:=\ffrac1{\sqrt{\alpha(P)}}x$$
be the corresponding $g_1$ spacelike or timelike unit vector. Similarly, if $\{e_1,e_2\}$ is an oriented $g_2$ orthonormal basis
for $\pi$, then 
$$\{\ffrac1{\sqrt\alpha}e_1,\ffrac1{\sqrt\alpha}e_2\}$$ is the corresponding oriented $g_1$ orthonormal basis for $\pi$. We then
have
$$J_{W_{g_1}}(\tilde x)=\ffrac1{\alpha(P)}J_{W_{g_2}}(x)\quad\text{and}\quad W_{g_1}(\pi)=\ffrac1{\alpha(P)}W_{g_2}(\pi)\,.
$$
The Lemma now follows as the Jordan normal forms are simply rescaled. \end{proof}

Theorem \ref{thm-1.1} shows that the notions we are studying are well defined in conformal geometry and justifies the notation we have
employed. Here is a brief guide to the remainder of the paper. In Section
\ref{Sect-2}, we will present some results concerning conformally spacelike and timelike Jordan Osserman manifolds. In Section \ref{Sect-3},
we will present some results concerning conformally spacelike and timelike Jordan Ivanov-Petrova manifolds. We conclude in Section
\ref{Sect-4} with some examples.

\section{Conformally Jordan Osserman manifolds}\label{Sect-2}

We begin with the following observation:

\begin{theorem}\label{thm-2.1}
If $(M,g)$ is Einstein, then $(M,g)$ is conformally spacelike (resp. timelike) Jordan Osserman if and only if
$(M,g)$ is pointwise spacelike (resp. timelike) Jordan Osserman.
\end{theorem}

\begin{proof}
If $(M,g)$ is Einstein, then Equation (\ref{eqn-1.b}) implies
$$g(J_W(x)y,z)=g(J_R(x)y,z)+\lambda\{g(y,z)g(x,x)-g(y,x)g(z,x)\}$$
where $\lambda$ is a suitably chosen constant. Thus
$$J_W(x)y=\left\{\begin{array}{lll}
0&\text{if}&y=x,\\
\{J_R(x)+\lambda g(x,x)\id\}y&\text{if}&y\perp x\,.
\end{array}\right.$$
Thus apart from the trivial eigenvalue $0$, the Jordan normal form of $J_W(x)$ and $J_R(x)$ are simply shifted by adding a scalar multiple of
the identity if $x$ is not a null vector. Theorem \ref{thm-2.1} is now immediate.\end{proof}

The classification is complete in certain settings:

\begin{theorem}\label{thm-2.2} Assume either that $(M,g)$ is an odd dimensional Riemannian manifold or that $(M,g)$ is
a Lorentzian manifold. Then
$(M,g)$ is conformally spacelike Jordan Osserman if and only if $(M,g)$ is conformally flat.
\end{theorem}

\begin{proof} We say $\mathcal{V}=(V,g,A)$ has {\it constant sectional curvature $\lambda$} if $A=\lambda A_{\id}$, i.e.
$$A(x,y,z,w)=\lambda\{g(x,w)g(y,z)-g(x,z)g(y,w)\}\,.$$
If $\mathcal{V}$ is Riemannian spacelike Jordan Osserman and if the dimension
$m$ is odd, then work of Chi \cite{Ch88} shows that $\mathcal{V}$ has constant sectional curvature. If $\mathcal{V}$ is
Lorentzian and spacelike Jordan Osserman, then results of Bla\v zi\'c, Bokan and Gilkey
\cite{BBG97} and of Garc\'{\i}a--R\'{\i}o, Kupeli and V\'azquez-Abal \cite{GKV97} shows that
$\mathcal{V}$ has constant sectional curvature.

If $A$ has constant sectional curvature $\lambda$ and if $x$ is not null, then
Equation (\ref{eqn-1.c}) shows
$$J_A(x)(y)=\left\{\begin{array}{lll}
0&\text{if }&y=x,\\
\lambda g(x,x)y&\text{if}&y\perp x\,.\end{array}\right.$$
Consequently $\Tr(J_A(x))=(m-1)\lambda g(x,x)$. Therefore, if $A\in\mathcal{W}(V,g)$, then necessarily $\lambda=0$ by Equation
(\ref{eqn-1.d}). Theorem
\ref{thm-2.2} now follows by applying these observations to $\mathcal{V}:=(T_PM,g_P,W_P)$.\end{proof}

One says that a manifold is spacelike (resp. timelike) Osserman if the eigenvalues of the Jacobi operator are constant on
the pseudo-sphere $S^+(T_PM)$ (resp. $S^-(T_PM)$ for any point $P\in M$. Theorem \ref{thm-2.2} extends to show that any Lorentzian
manifold
$(M,g)$ which is spacelike (resp. timelike) Osserman is conformally flat.

Any local rank $1$ Riemannian symmetric space is necessarily conformally Osserman since the group of local isometries acts transitively on
the unit sphere bundle. We conjecture that the converse holds; this is the analogue of the Osserman conjecture in this setting:

\begin{conjecture}\label{conj-2.3}
A connected Riemannian manifold $(M,g)$ is conformally Osserman if and only if $(M,g)$ is locally conformally equivalent to a rank $1$
symmetric space.
\end{conjecture}

We shall see in Section \ref{Sect-4} that this conjecture fails in the higher signature setting.

\section{Conformally Jordan Ivanov-Petrova manifolds}\label{Sect-3}

The classification is almost complete in the Riemannian setting:

\begin{theorem}\label{thm-3.1} 
Let $(M,g)$ be a conformally spacelike Jordan Ivanov-Petrova Riemannian manifold of dimension $m\ne 3,7$. Then $(M,g)$ is
conformally flat.
\end{theorem}

\begin{proof} Suppose first $m\ge5$ and $m\ne 7$. We apply results of \cite{Gi99,GLS99} to see that any Riemannian
Ivanov-Petrova algebraic curvature tensor in these dimensions has rank $2$. Such tensors are classified. Let $P\in M$. There
exists a self-adjoint isometry
$\phi_P$ of
$T_PM$ with $\phi_P^2=\id$ so that $W_P=\lambda R_{\phi_P}$ where $R_{\phi_P}$
is given by Equation (\ref{eqn-1.a}). We may then use Equation (\ref{eqn-1.c})
to see:
\begin{equation}\label{eqn-3.a}
J_W(x)y=\lambda g_P(\phi_Px,x)\phi_Py\quad\text{if}\quad y\perp\phi_Px\,.
\end{equation}
Decompose $T_PM=T_P^+M\oplus T_P^-M$ into the $\pm1$ eigenspaces of $\phi_P$.
Let $e^\pm$ be unit vectors in
$T_P^\pm M$. Set $a^\pm:=\dim T_P^\pm M$. Then Equation (\ref{eqn-3.a}) implies
that
\begin{equation}\label{eqn-3.b}
\Tr\{J_W(e^+)\}=\lambda(a^+-1-a^-)\quad\text{and}\quad\Tr\{J_W(e^-)\}=\lambda(a^--1-a^+)\,.
\end{equation}

If $\phi=\pm\id$, then $W_P$ has constant sectional curvature and the argument given to establish Theorem \ref{thm-2.2} shows $W_P=0$. Thus
we may assume that
$a^+\ge1$ and $a^-\ge1$. By Equation (\ref{eqn-1.d}),
$\Tr\{J_W(x)\}=0$ for any $x$. Thus we have 
$$(a^+-a^--1)\lambda =0\quad\text{and}\quad(a^--a^+-1)\lambda =0\,.$$
 Adding these two equations implies
$-2\lambda=0$ and hence $W_P=0$. This establishes the Lemma except when $m=4$.

We complete the proof of the Lemma by dealing with the exceptional case $m=4$. We follow the discussion in Ivanov-Petrova \cite{IvPe98} to
see that either $W$ has the form given in Equation (\ref{eqn-3.a}), in which case the argument given above shows $W_P=0$, or that there
exists an orthonormal basis
$\{e_1,e_2,e_3,e_4\}$ for
$T_PM$ so that the non-zero components of $W$ are given by:
\begin{equation}\label{eqn-3.c}
\begin{array}{llll}
 W_{1212}=a_1,\hfill& W_{1234}=\phantom{-}a_2,\hfill&
 W_{1313}=a_2,\hfill& W_{1324}=-a_1,\hfill\\
 W_{1414}=a_2,\hfill& W_{1423}=\phantom{-}a_1,\hfill&
 W_{2323}=a_2,\hfill& W_{2314}=\phantom{-}a_1,\hfill\\\
 W_{2424}=a_2,\hfill& W_{2413}=-a_1,\hfill&
 W_{3434}=a_1,\hfill&W_{3412}=\phantom{-}a_2\,.
\end{array}
\end{equation}
where $a_2+2a_1=0$. Since $\rho_W(e_1,e_1)=-2a_2-a_1=0$ by Equation
(\ref{eqn-1.d}), we conclude $a_1=a_2=0$, which once again implies
$W_P=0$.
\end{proof}

There are analogous results in the higher signature setting, although with slightly more restrictive hypotheses.

\begin{theorem}\label{thm-3.2}
Let $(M,g)$ be a connected pseudo-Riemannian manifold of signature $(p,q)$ which is conformally spacelike Ivanov-Petrova. Assume that
$q\ge11$, that
$p\le\frac{q-6}4$, and that $\{q,q+1,...,q+p\}$ does not contain a power of $2$. Then either $W(\pi)$ is nilpotent for every spacelike $2$
plane or $(M,g)$ is conformally flat.
\end{theorem}

\begin{proof} Results of \cite{GiZa02a,GiZa02b,St03} show that there exists a normalizing constant $\lambda$ so that $W_P=\lambda W_\phi$
where $W_\phi$ is
given by Equation (\ref{eqn-1.a}) where one of the following conditions holds:
\begin{enumerate}
\item $\phi^2=\id$ and $\phi$ is a self-adjoint isometry of $T_PM$.
\item $\phi^2=-\id$ and $\phi$ is a self-adjoint para-isometry of $T_PM$.
\item $\phi^2=0$.
\end{enumerate}
If $\phi^2=0$, then $W(\pi)$ is always nilpotent. We complete the proof by showing that either (1) or (2) imply $\lambda=0$.

Suppose $\phi$ is a self-adjoint isometry of $T_PM$ with $\phi^2=\id$. As in the proof of
Theorem \ref{thm-3.1}, we decompose 
$T_PM=T_P^+M\oplus T_P^-M$ into the $\pm1$
eigenspaces of $\phi$.  Again, we set $a^\pm=\dim T_P^\pm$ were we may suppose $a^+\ge1$ and $a^-\ge1$. These
eigenspaces are orthogonal with respect to the metric $g_P$ and thus the restriction of the metric to each
eigenspace is non-degenerate. Thus we may choose vectors
$e^\pm\in T_P^\pm M$ so $g_P(e^\pm,e^\pm)=\varepsilon^\pm\ne0$. Equation (\ref{eqn-3.b}) then extends to become
$$
\Tr\{J_W(e^+)\}=\varepsilon^+\lambda(a^+-1-a^-)\quad\text{and}\quad\Tr\{J_W(e^-)\}=\varepsilon^-\lambda(a^--1-a^+)\,.
$$
We argue as in the proof of Theorem \ref{thm-3.1} to see that this implies $\lambda=0$. 

If $\phi$ is a para-isometry, we complexify.
Replacing $\phi$ by $\tilde\phi:=\sqrt{-1}\phi$ and applying the argument given above to the $g$ self-adjoint (complex) isometry
$\tilde\phi$ to see that
$\sqrt{-1}\lambda=0$ and thus, again,
$W_P=0$.
\end{proof}

This result, together with the examples in the subsequent section, motivates the following:

\begin{conjecture}\label{conj-3.3}
 Let $(M,g)$ be a conformally spacelike Ivanov-Petrova manifold. If $(M,g)$ is not
conformally flat, then $W(\pi)$ is nilpotent for any oriented spacelike $2$ plane $\pi$.
\end{conjecture}

\section{Examples}\label{Sect-4}

Theorem \ref{thm-2.2} shows that Conjecture \ref{conj-2.3} holds if $m$ is odd. The situation is considerably more complicated in the higher
signature setting. The following family of manifolds \cite{GIZ03} is useful in this setting. It also shows there are conformally spacelike
Jordan Ivanov-Petrova manifolds which are not conformally flat. Let
$p\ge2$. Introduce coordinates
$(x_1,...,x_p,y_1,...,y_p)$ on $\mathbb{R}^{2p}$ and let $f=f(x_1,...,x_p)$ be a smooth function on $\mathbb{R}^{2p}$. Define a neutral
signature metric $g_f$ on $\mathbb{R}^{2p}$ by setting
$$g_f(\partial_i^x,\partial_j^x):=\partial_i^xf\cdot\partial_j^xf,\quad
 g_f(\partial_i^x,\partial_j^y)=g_f(\partial_j^y,\partial_i^x)=\delta_{ij},\quad
 g_f(\partial_i^y,\partial_j^y)=0\,.$$
Let $H=(H_{ij})\in M_p(\mathbb{R})$ be the Hessian where $H_{ij}=\partial_i^x\partial_j^x$.

\begin{theorem}\label{thm-4.1}
Let $(M,g_f)$ be as defined above. Assume that $p\ge3$.
\begin{enumerate}
\item Assume that $H$ is definite. If $x$ is not null, then $J_W(x)$ has rank $p-1$
and $J_W(x)^2=0$. Thus $(M,g_f)$ is conformally spacelike and timelike Jordan Osserman.
\item If $H$ is indefinite, then $(M,g_f)$ is neither conformally spacelike Jordan Osserman nor conformally timelike Jordan Osserman.
\item Assume that $H$ is non-degenerate. If $\pi$ is an oriented spacelike or timelike $2$ plane, then $\Rank(W(\pi))=2$ and $W(\pi)^2=0$.
Thus
$(M,g_f)$ is conformally spacelike and timelike Jordan Ivanov-Petrova.
\end{enumerate}
\end{theorem}

\begin{proof} We showed in \cite{GIZ03} that $(M,g_f)$ was Ricci flat. Consequently, $W=R$. The assertions of Theorem \ref{thm-4.1} now
follow from the corresponding assertions for $J_R(x)$ and for $R(\pi)$ which were established in \cite{GIZ03}. \end{proof}

The manifolds of Theorem \ref{thm-4.1} have a Jacobi operator and skew-symmetric curvature operator which are nilpotent of order $2$. There
are also manifolds where the Jacobi operator and skew-symmetric curvature operator are nilpotent of order $3$. Let
$(u_1,...,u_s,t_1,...,t_s,w_1,...,w_s)$ be coordinates on
$\mathbb{R}^{3s}$ for
$s\ge2$. Let $f_i(x)$ be smooth functions on $\mathbb{R}$ and set $F(u_1,...,u_s):=f_1(u_1)+...+f_s(u_s)$. Define a metric of signature
$(2s,s)$ on $M_F:=\mathbb{R}^{3s}$ by setting
$$\begin{array}{ll}
g_F(\partial_i^u,\partial_j^u)=-2\delta_{ij}F(u)-2\delta_{ij}\sum_k u_kt_k,&
g_F(\partial_i^u,\partial_j^v)=g_F(\partial_j^v,\partial_i^u)=\delta_{ij},\\
g_F(\partial_i^u,\partial_j^t)=g_F(\partial_j^t,\partial_i^u)=0,&
g_F(\partial_i^t,\partial_j^t)=-\delta_{ij},\vphantom{\vrule height 11pt}\\
g_F(\partial_i^t,\partial_j^v)=g_F(\partial_j^v,\partial_i^t)=0,&
g_F(\partial_i^v,\partial_j^v)=0\,.\vphantom{\vrule height 11pt}
\end{array}$$

\begin{theorem}\label{thm-4.2}
Let $(M_F,g_F)$ be as defined above where $s\ge2$. 
\begin{enumerate}\item Let $x$ be spacelike. Then $J_W(x)$ has rank $2s-2$, $J_W(x)^2$ has rank $s-1$, and\newline
$J_W(x)^3=0$. Consequently $(M_F,g_F)$ is conformally spacelike Jordan Osserman. However, $(M_F,g_F)$ is not conformally timelike Jordan
Osserman.
\item Let $\pi$ be an oriented spacelike $2$ plane. Then $W(\pi)$ has rank $4$, $W(\pi)^2$ has rank $2$, and $W(\pi)^3=0$.
Consequently $(M_F,g_F)$ is conformally spacelike Jordan Ivanov-Petrova. However, $(M_F,g_F)$ is not conformally timelike Jordan
Ivanov-Petrova.
\end{enumerate}
\end{theorem}

\begin{proof} We showed in \cite{GN03} that $(M_F,g_F)$ is Ricci flat. Consequently, $R=W$. The
assertions of the Theorem now follow from the corresponding assertions for $J_R$ which were established in \cite{GN03} and for $W(\pi)$
which were established in
\cite{GNSV03}.\end{proof}

\section*{Acknowledgments} Research of N. Bla\v zi\'c partially supported by the DAAD (Germany) and MNTS Project \#1854 (Srbija). Research of
P. Gilkey partially supported by the MPI (Leipzig). Research of S. Nik\v cevi\'c partially supported by the DAAD (Germany) and MMTS \#1646
(Srbija). Research of U. Simon partially supported by DFG-Si163/03. The authors wish to express their thanks to the Technische
Universit\"at Berlin where much of the research reported here was conducted.

\end{document}